\documentclass[draft,11pt,ukrainian]{article}

\usepackage{amsfonts,amsmath,amsthm,amscd,amssymb,latexsym,cite,verbatim,texdraw,floatflt,caption2,pb-diagram}

    \usepackage[T2A]{fontenc} 
    \usepackage[cp1251]{inputenc}
    \usepackage[russian]{babel}

\theoremstyle{plain}
\newtheorem{theorem}{Теорема}[section]
\newtheorem{lemma}{Лема}[section]
\newtheorem{proposition}{Твердження}[section]

\theoremstyle{definition}
\newtheorem{example}{Приклад}[section]
\newtheorem{remark}{Зауваження}[section]
\newcommand{\keywords}{\textbf{Ключові слова. }}
\newcommand{\subjclass}{\textbf{MSC 2010. }}
\renewcommand{\abstract}{\textbf{Анотація. }}

\,
\,

\numberwithin{equation}{section}

\begin{document}

\title{Про моногенні функції, визначені в різних комутативних алгебрах
}

\author{Віталій С. Шпаківський}



\date{}

\maketitle

\begin{abstract}
Встановлено відповідність між моногенною функцією в довільній
скінченновимірній комутативній асоціативній алгебрі і скінченним
набором  моногенних функцій в спеціальній комутативній асоціативній
алгебрі.
\end{abstract}



\section{Вступ}

Напевно першим хто використав аналітичні функції, що приймають
значення в комутативній алгебрі для побудови розв'язків тривимірного
рівняння Лапласа був П.~Кетчум \cite{Ketchum-28}. Він показав, що
кожна аналітична функція $\Phi(\zeta)$ змінної
$\zeta=xe_1+ye_2+ze_3$ задовольняє тривимірне рівняння Лапласа, якщо
лінійно незалежні елементи $e_1, e_2, e_3$ комутативної алгебри
задовольняють умову
\begin{equation}\label{garmonichnyj_bazys-ogljad}
     e_1^2+e_2^2+e_3^2=0\,,
\end{equation}
оскільки
\begin{equation}\label{garm}
\Delta_3\Phi:=\frac{{\partial}^{2}\Phi}{{\partial x}^{2}}+
\frac{{\partial}^{2}\Phi}{{\partial y}^{2}}+
\frac{{\partial}^{2}\Phi}{{\partial z}^{2}}\equiv{\Phi}''(\zeta) \
(e_1^2+e_2^2+e_3^2)=0\,,\medskip
\end{equation}
де $\Phi'':=(\Phi')'$ і $\Phi'(\zeta)$ визначається рівністю
 $d\Phi=\Phi'(\zeta)d\zeta$.

Узагальнюючи П. Кетчума,  М. Рошкулець \cite{Rosculet,Rosculet-56}
використовував аналітичні функції зі значеннями в комутативних
алгебрах для дослідження рівнянь вигляду
\begin{equation}\label{dopolnenije----1}
\mathcal{L}_NU(x,y,z):=\sum\limits_{\alpha+\beta+\gamma=N}C_{\alpha,\beta,\gamma}\,
\frac{\partial^N U} {\partial x^\alpha\,\partial y^\beta\,\partial
z^\gamma}=0, \quad C_{\alpha,\beta,\gamma}\in\mathbb{R}.
\end{equation}

Розглядаючи змінну $\zeta=xe_1+ye_2+ze_3$ і аналітичну функцію
$\Phi(\zeta)$, отримуємо наступну рівність для мішаної похідної:
\begin{equation}\label{garm+1}
\frac{\partial^{\alpha+\beta+\gamma}\Phi} {\partial
x^\alpha\,\partial y^\beta\,\partial z^\gamma}= e_1^\alpha\,
e_2^\beta\,
e_3^\gamma\,\Phi^{(\alpha+\beta+\gamma)}(\zeta)=e_1^\alpha\,
e_2^\beta\, e_3^\gamma\,\Phi^{(N)}(\zeta).
\end{equation}
Підставляючи (\ref{garm+1}) в рівняння (\ref{dopolnenije----1}),
маємо рівність
$$
\mathcal{L}_N\Phi(\zeta)=\Phi^{(N)}(\zeta)
\sum\limits_{\alpha+\beta+\gamma=N}C_{\alpha,\beta,\gamma}\,e_1^\alpha\,e_2^\beta
\,e_3^\gamma\,.
$$

Приходимо до висновку, що для виконання рівності
$\mathcal{L}_N\Phi(\zeta)=0$ елементи алгебри $e_1=1,e_2,e_3$ мають
задовольняти \textit{характеристичне рівняння}
\begin{equation}\label{dopolnenije----2}
\mathcal{X}(1,e_2,e_3):=\sum\limits_{\alpha+\beta+\gamma=N}C_{\alpha,\beta,\gamma}\,e_2^\beta\,
e_3^\gamma=0\,.
\end{equation}

Якщо ліву частину рівняння (\ref{dopolnenije----2}) розкласти за
базисом алгебри, то характеристичне рівняння
(\ref{dopolnenije----2}) рівносильне \textit{характеристичній
системі рівнянь}, породженій рівнянням (\ref{dopolnenije----2}).

Таким чином, при виконанні умови (\ref{dopolnenije----2}) кожна
аналітична функція $\Phi$  зі значеннями в довільній комутативній
асоціативній алгебрі задовольняє рівняння (\ref{dopolnenije----1}),
і, відповідно, усі дійснозначні компоненти функції  $\Phi$ є
розв'язками рівняння (\ref{dopolnenije----1}).

В роботі \cite{Pogor-Ramon-Shap} розглядаються диференціальні
рівняння в частинних похідних від декількох змінних і наведено ряд
прикладів на застосування описаного вище методу.

І. Мельниченко \cite{Mel'nichenko75} запропонував розглядати в
рівностях \eqref{garm} і \eqref{garm+1} функції $\Phi$, двічі
диференційовні за Гато, при цьому описав усі базиси $\{e_1, e_2,
e_3\}$ тривимірних комутативних алгебр з одиницею над полем
$\mathbb{C}$, які задовольняють рівність
(\ref{garmonichnyj_bazys-ogljad}), див. \cite{Plaksa}.

Для цих тривимірних комутативних алгебр, асоційованих з тривимірним
рівнянням Лапласа, в роботах \cite{Pl-Shp1,Pl-Pukh,Pukh} отримано
конструктивний опис усіх моногенних (тобто неперервних і
диференційовних за Гато) функцій за допомогою трьох відповідних
голоморфних функцій комплексної змінної.

В роботах \cite{Pl-Shp-Algeria,Pl-Pukh-Analele} встановлено
конструктивний опис моногенних функцій (зв'язаних з рівнянням
$\Delta_3\Phi=0$) зі значеннями в деяких $n$-вимірних комутативних
алгебрах за допомогою відповідних $n$ голоморфних функцій
комплексної змінної і, спираючись на одержані представлення
моногенних функцій, доведено аналоги ряду класичних результатів
комплексного аналізу.

Нарешті в роботі \cite{Shpakivskyi-2014} отримано конструктивний
опис моногенних функцій (зв'язаних з рівнянням
(\ref{dopolnenije----1})) зі значеннями в довільній комутативній
асоціативній алгебрі над полем $\mathbb{C}$ за допомогою голоморфних
функцій комплексної змінної.

У цій роботі буде показано, що для побудови розв'язків рівняння
(\ref{dopolnenije----1}) у вигляді компонент моногенних функцій зі
значеннями в скінченновимірних комутативних асоціативних алгебрах
достатньо обмежитись вивченням моногенних функцій у алгебрах певного
виду.

\section{Алгебра $\mathbb{A}_n^m$}

Нехай $\mathbb{N}$ --- множина натуральних
 чисел і  $m,n\in\mathbb{N}$ такі, що $m\leq n$.
Нехай $\mathbb{A}_n^m$ --- довільна комутативна асоціативна  алгебра
з одиницею над полем комплексних чисел $\mathbb{C}$.
 Е.~Картан \cite[с.~33]{Cartan} довів, що в алгебрі $\mathbb{A}_n^m$
 існує базис $\{I_k\}_{k=1}^{n}$, який задовольняє наступні правила множення:

\vskip3mm 1.  \,\,  $\forall\, r,s\in[1,m]\cap\mathbb{N}\,:$ \qquad
$I_rI_s=\left\{
\begin{array}{rcl}
0 &\mbox{при} & r\neq s,\vspace*{2mm} \\
I_r &\mbox{при} & r=s;\\
\end{array}
\right.$

\vskip5mm

2. \,\,  $\forall\, r,s\in[m+1,n]\cap\mathbb{N}\,:$ \qquad $I_rI_s=
\sum\limits_{k=\max\{r,s\}+1}^n\Upsilon_{r,k}^{s}I_k$\,;

\vskip5mm

3.\,\, $\forall\, s\in[m+1,n]\cap\mathbb{N}$  $\exists!\;
 u_s\in[1,m]\cap\mathbb{N}$ \;$\forall\,
 r\in[1,m]\cap\mathbb{N}\,:$\;\;

$$I_rI_s=\left\{
\begin{array}{ccl}
0 \;\;\mbox{при}\;\;  r\neq u_s\,,\vspace*{2mm}\\
I_s\;\;\mbox{при}\;\;  r= u_s\,. \\
\end{array}
\right.\medskip
$$
Крім того, структурні константи $\Upsilon_{r,k}^{s}\in\mathbb{C}$
задовольняють умови асо\-ціа\-тив\-ності: \vskip2mm (A\,1).\quad
$(I_rI_s)I_p=I_r(I_sI_p)$ \;
$\forall\,r,s,p\in[m+1,n]\cap\mathbb{N}$; \vskip2mm (A\,2).\quad
$(I_uI_s)I_p=I_u(I_sI_p)$ \; $\forall\, u\in[1,m]\cap\mathbb{N}$\;
 $\forall\, s,p\in[m+1,n]\cap\mathbb{N}$. \vskip2mm Очевидно, що перші $m$
базисних векторів $\{I_u\}_{u=1}^m$ є ідемпотентами і породжують
напівпросту підалгебру $S$ алгебри $\mathbb{A}_n^m$, а вектори
$\{I_r\}_{r=m+1}^n$ породжують нільпотентну підалгебру $N$ цієї
алгебри. З правил множення алгебри $\mathbb{A}_n^m$ випливає, що
$\mathbb{A}_n^m$ є напівпрямою сумою $m$-вимірної напівпростої
підалгебри $S$ і $(n-m)$-вимірної нільпотентної підалгебри $N$,
тобто
$$\mathbb{A}_n^m=S\oplus_s N.
$$
Одиницею алгебри $\mathbb{A}_n^m$ є елемент $1=\sum_{u=1}^mI_u$.

Алгебра $\mathbb{A}_n^m$ містить $m$ максимальних ідеалів
$$\mathcal{I}_u:=\Biggr\{\sum\limits_{k=1,\,k\neq u}^n\lambda_kI_k:\lambda_k\in
\mathbb{C}\Biggr\}, \quad  u=1,2,\ldots,m,
$$
перетином яких є радикал
\begin{equation}\label{radyk}
\mathcal{R}:= \Bigr\{\sum\limits_{k=m+1}^n\lambda_kI_k:\lambda_k\in
\mathbb{C}\Bigr\}.
\end{equation}

Визначимо $m$ лінійних функціоналів
$f_u:\mathbb{A}_n^m\rightarrow\mathbb{C}$ рівностями
\begin{equation}\label{funkc}
f_u(I_u)=1,\quad f_u(\omega)=0 \,\,\,\forall\, \omega\in\mathcal{I}_u,\quad u=1,2,\ldots,m.
\end{equation}
Оскільки ядрами функціоналів $f_u$ є відповідно максимальні ідеали
 $\mathcal{I}_u$, то ці функціонали є також
 неперервними і мультиплікативними
  (див. \cite[с. 147]{Hil_Filips}).

\section{Моногені функції}

Нехай
\begin{equation}\label{e_1_e_2_e_3+}
e_1=1,\quad e_2=\sum\limits_{r=1}^na_rI_r,\quad
e_3=\sum\limits_{r=1}^nb_rI_r
\end{equation}
при $a_r,b_r\in\mathbb{C}$ --- трійка векторів в алгебрі
$\mathbb{A}_n^m$, які лінійно незалежні над полем $\mathbb{R}$. Це
означає, що рівність
$$\alpha_1e_1+\alpha_2e_2+\alpha_3e_3=0,\quad \alpha_1,\alpha_2,
\alpha_3\in\mathbb{R},$$ виконується тоді і тільки тоді, коли
$\alpha_1=\alpha_2= \alpha_3=0$.

Нехай $\zeta:=xe_1+ye_2+ze_3$, де $x,y,z\in\mathbb{R}$. Очевидно, що
$\xi_u:=f_u(\zeta)=x+ya_u+zb_u$,\, $u=1,2,\ldots,m$.
 Виділимо в алгебрі $\mathbb{A}_n^m$ лінійну оболонку
 $E_3:=\{\zeta=xe_1+ye_2+ze_3:\,\, x,y,z\in\mathbb{R}\}$, породжену
 векторами $e_1,e_2,e_3$.

 Далі істотним є припущення:
 $f_u(E_3)=\mathbb{C}$ при всіх $u=1,2,\ldots,m$, де $f_u(E_3)$ --- образ множини $E_3$ при відображенні $f_u$.
Очевидно, що це має місце тоді і тільки тоді, коли при кожному
фіксованому $u=1,2, \ldots, m$ хоча б одне з чисел $a_u$ чи $b_u$
належить
 $\mathbb{C}\setminus\mathbb{R}$. В теоремі 7.1  роботи \cite{Shpakivskyi-2014}
встановлено підклас рівнянь вигляду (\ref{dopolnenije----1}) для
яких умова $f_u(E_3)=\mathbb{C}$ виконується при всіх
$u=1,2,\ldots,m$.

 Області $\Omega$
тривимірного простору $\mathbb{R}^{3}$ поставимо у відповідність
область $\Omega_{\zeta}:=\{\zeta =xe_{1}+ye_{2}+ze_{3}: (x, y, z)\in
\Omega\}$ в $E_{3}$.

Неперервну функцію $\Phi:\Omega_{\zeta}\rightarrow\mathbb{A}_n^m$
називатимемо {\it моногенною}\/ в області $\Omega_{\zeta}\subset
E_{3}$, якщо $\Phi$ диференційовна за Гато в кожній точці цієї
області, тобто якщо для кожного $\zeta\in\Omega_{\zeta}$ існує
елемент $\Phi'(\zeta)$ алгебри $\mathbb{A}_n^m$ такий, що
виконується рівність
$$\lim\limits_{\varepsilon\rightarrow 0+0}
\left(\Phi(\zeta+\varepsilon h)-\Phi(\zeta)\right)\varepsilon^{-1}=
h\Phi'(\zeta)\quad\forall\, h\in E_{3}.
$$
$\Phi'(\zeta)$ називається {\it похідною Гато}\/ функції $\Phi$ в
точці $\zeta$.

Розглянемо розклад функції
$\Phi:\Omega_{\zeta}\rightarrow\mathbb{A}_n^m$ за базисом
$\{I_k\}_{k=1}^n$:
\begin{equation}\label{rozklad-Phi-v-bazysi}
\Phi(\zeta)=\sum_{k=1}^n U_k(x,y,z)\,I_k\,.
 \end{equation}

У випадку, коли функції $U_k:\Omega\rightarrow\mathbb{C}$ є
$\mathbb{R}$-диференційовними в області $\Omega$, тобто для
довільного $(x,y,z)\in\Omega$
$$U_k(x+\Delta x,y+\Delta y,z+\Delta z)-U_k(x,y,z)=
\frac{\partial U_k}{\partial x}\,\Delta x+ \frac{\partial
U_k}{\partial y}\,\Delta y+\frac{\partial U_k}{\partial z}\,\Delta
z+$$ $$+\,o\left(\sqrt{(\Delta x)^2+(\Delta y)^2+(\Delta
z)^2}\,\right), \qquad (\Delta x)^2+(\Delta y)^2+(\Delta z)^2\to
0\,,$$ функція $\Phi$ моногенна в області $\Omega_{\zeta}$ тоді і
тільки тоді, коли у кожній точці області
 $\Omega_{\zeta}$ виконуються умови:
\begin{equation}\label{Umovy_K-R}
\frac{\partial \Phi}{\partial y}=\frac{\partial \Phi}{\partial
x}\,e_{2}\,,\quad \frac{\partial \Phi}{\partial z}=\frac{\partial
\Phi}{\partial x}\,e_{3}\,.
\end{equation}

Відмітимо, що розклад резольвенти має вигляд
\begin{equation}\label{rozkl-rezol-A_n^m}
(te_1-\zeta)^{-1}=\sum\limits_{u=1}^m\frac{1}{t-\xi_u}\,I_u+
 \sum\limits_{s=m+1}^{n}\sum\limits_{k=2}^{s-m+1}\frac{Q_{k,s}}
 {\left(t-\xi_{u_{s}}\right)^k}\,I_{s}\,
  \end{equation}
  $$ \forall\,t\in\mathbb{C}:\,
t\neq \xi_u,\quad u=1,2,\ldots,m,$$ де $Q_{k,s}$ визначені
наступними рекурентними співвідношеннями:
$$
Q_{2,s}:=T_{s}\,,\quad
Q_{k,s}=\sum\limits_{r=k+m-2}^{s-1}Q_{k-1,r}\,B_{r,\,s}\,,\;
\;\;k=3,4,\ldots,s-m+1.
$$
при
$$T_s:=ya_s+zb_s\,,\;\;B_{r,s}:=\sum\limits_{k=m+1}^{s-1}T_k
\Upsilon_{r,s}^k\,,\; \;\;s=m+2,\ldots,n,$$
 а натуральні числа $u_s$  визначені у правилі 3 таблиці множення алгебри $\mathbb{A}_n^m$.

Із співвідношень  (\ref{rozkl-rezol-A_n^m}) випливає, що точки
 $(x,y,z)\in\mathbb{R}^3$, які відповідають необоротним елементам
 $\zeta\in\mathbb{A}_n^m$, лежать на прямих
 \begin{equation}\label{L-u}
  L_u:\quad\left\{
\begin{array}{r}x+y\,{\rm Re}\,a_u+z\,{\rm Re}\,b_u=0,\vspace*{2mm} \\
y\,{\rm Im}\,a_u+z\,{\rm Im}\,b_u=0 \\
\end{array} \right.
 \end{equation}
в тривимірному просторі $\mathbb{R}^3$.

Нехай область $\Omega\subset \mathbb{R}^{3}$ є опуклою в напрямку
прямих $L_u$,\, $u=1,2,\ldots, m$. Позначимо через $D_u$ область
комплексної площини $\mathbb{C}$\,, на яку область $\Omega_\zeta$
відображається функціоналом $f_u$.

\textbf{Теорема А}\cite{Shpakivskyi-2014}. \textit{Нехай область
$\Omega\subset \mathbb{R}^{3}$ є опуклою в напрямку прямих $L_u$ і
 $f_u(E_3)=\mathbb{C}$ при всіх
 $u=1,2,\ldots, m$.
Тоді кожна моногенна функція
$\Phi:\Omega_{\zeta}\rightarrow\mathbb{A}_n^m$ подається у вигляді
 \begin{equation}\label{Teor--1}
\Phi(\zeta)=\sum\limits_{u=1}^mI_u\,\frac{1}{2\pi
i}\int\limits_{\Gamma_u} F_u(t)(t-\zeta)^{-1}\,dt+
\sum\limits_{s=m+1}^nI_s\,\frac{1}{2\pi i}\int\limits_
{\Gamma_{u_s}}G_s(t)(t-\zeta)^{-1}\,dt,
 \end{equation}
де $F_u$ --- деяка голоморфна функція в області $D_u$ і $G_s$ ---
 деяка голоморфна функція
в області $D_{u_s}$, а $\Gamma_q$ --- замкнена жорданова спрямлювана
крива, яка лежить в області $D_q$, охоплює  точку $\xi_q$ і не
містить точок $\xi_{\ell}$\,, $\ell=1,2,\ldots, m$,\,$\ell\neq q$.}

Оскільки за умов теореми \textbf{А} кожна моногенна функція
$\Phi:\Omega_{\zeta}\rightarrow \mathbb{A}_n^m$ продовжується до
функції, моногенної в області
 \begin{equation}\label{Pi}
\Pi_\zeta:=\{\zeta\in E_3:f_u(\zeta)=D_u\,,\,u=1,2,\ldots,m\},
 \end{equation}
то надалі будемо розглядати моногенні функції $\Phi$, визначені в
областях виду $\Pi_\zeta$\,.

\section{Характеристичне рівняння в різних комутативних алгебрах}

Скажемо, що система поліноміальних над полем $\mathbb{C}$ рівнянь
$Q_1$ \textit{редукується} до системи поліноміальних рівнянь $Q_2$,
якщо система $Q_2$ отримується з системи $Q_1$ шляхом відкидання
деякої кількості рівнянь. В свою чергу, система $Q_2$ \textit{є
редукцією системи} $Q_1$. Відмітимо, що для заданої системи
поліноміальних рівнянь $Q_1$ редукована система $Q_2$ не єдина.
Очевидним є наступне твердження.

\begin{proposition}\label{Prop1}
Нехай система
 поліноміальних рівнянь $Q_1$ з комплексними невідомими $t_1,t_2,\ldots,t_n$ має розв'язки і
$Q_2$ --- будь-яка її редукована система з невідомими
$t_{i_1},t_{i_2},\ldots,t_{i_k}$, де $i_1,i_2,\ldots,i_k$, $k\leq
n$, --- попарно різні елементи множини $\{1,2,\ldots,n\}$. Тоді усі
$t_{i_1},t_{i_2},\ldots,t_{i_k}$, які задовольняють систему $Q_1$ є
розв'язками системи $Q_2$.
\end{proposition}

Наприклад, система рівнянь
\begin{equation}\label{s1}
\begin{array}{l}
1+a_1^2+b_1^2=0,\vspace*{2mm} \\
1+a_2^2+b_2^2=0,\vspace*{2mm} \\
a_2a_3+b_2b_3=0\\
\end{array}
 \end{equation}
редукується до системи рівнянь
\begin{equation}\label{s2}
\begin{array}{l}
1+a_2^2+b_2^2=0,\vspace*{2mm} \\
a_2a_3+b_2b_3=0.\\
\end{array}
 \end{equation}
Твердження \ref{Prop1} означає, що всі значення $a_2,b_2,a_3,b_3$,
які задовольняють систему (\ref{s1}) є розв'язками системи
(\ref{s2}).

 Встановимо допоміжні твердження.

\begin{lemma}\label{lem-1-umb}
Нехай в алгебрі $\mathbb{A}_n^m$ існує трійка лінійно незалежних над
$\mathbb{R}$ векторів $1,e_2,e_3$, які задовольняють характеристичне
рівняння (\ref{dopolnenije----2}). Тоді для кожного
$u\in\{1,2,\ldots,m\}$ характеристична система, породжена рівнянням
$\mathcal{X}(I_u,e_2I_u,e_3I_u)=0$ є редукцією характеристичної
системи, породженої рівнянням (\ref{dopolnenije----2}).
\end{lemma}

\begin{proof} Нехай ліва частина рівняння (\ref{dopolnenije----2}) в
базисі алгебри має вигляд
$$
\mathcal{X}(1,e_2,e_3)=\sum\limits_{\alpha+\beta+\gamma=N}C_{\alpha,\beta,\gamma}\,e_2^\beta\,
e_3^\gamma=\sum\limits_{k=1}^nV_k\,I_k=0.
$$
Відповідно, характеристична система, породжена рівнянням
(\ref{dopolnenije----2}), має вигляд
\begin{equation}\label{+11}
\begin{array}{l}
V_1=0,\vspace*{2mm} \\
\ldots\ldots\vspace*{2mm} \\
V_n=0.
\end{array}
 \end{equation}

Тепер розглянемо характеристичну систему, породжену рівнянням
$\mathcal{X}(I_u,e_2I_u,e_3I_u)=0$. Маємо
$$
\mathcal{X}(I_u,e_2I_u,e_3I_u)=\sum\limits_{\alpha+\beta+\gamma=N}C_{\alpha,\beta,\gamma}\,I_u\,(e_2I_u)^\beta\,
(e_3I_u)^\gamma=
$$
\begin{equation}\label{+1}
=I_u\sum\limits_{\alpha+\beta+\gamma=N}C_{\alpha,\beta,\gamma}\,e_2^\beta\,
e_3^\gamma=I_u\sum\limits_{k=1}^nV_k\,I_k=V_u+I_u\sum\limits_{k=m+1}^nV_kI_k=0.
 \end{equation}
Відповідно до правила 3 таблиці множення алгебри $\mathbb{A}_n^m$
добуток $I_u\sum\limits_{k=m+1}^nV_kI_k$ належить радикалу
$\mathcal{R}$. Таким чином, рівняння (\ref{+1}) рівносильне такій
характеристичній системі:
\begin{equation}\label{+111}
\begin{array}{l}
V_u=0,\vspace*{2mm} \\
\ldots\ldots\vspace*{2mm} \\
V_k=0 \quad  \forall\,k\in\{m+1,\ldots,n\}\, :\, I_uI_k=I_k\,.
\end{array}
 \end{equation}
Очевидно, що система (\ref{+111}) є редукцією системи (\ref{+11}).

\end{proof}

Позначимо через ${\rm Rad}\,e_2$ частину вектора $e_2$ з розкладу
(\ref{e_1_e_2_e_3+}), яка міститься в його радикалі, тобто ${\rm
Rad}\,e_2:=\sum\limits_{r=m+1}^na_rI_r\,$. Аналогічно, ${\rm
Rad}\,e_3:=\sum\limits_{r=m+1}^nb_rI_r\,$.

\begin{lemma}\label{lem-2-umb}
Нехай в алгебрі $\mathbb{A}_n^m=S\oplus_s N$ існує трійка лінійно
незалежних над $\mathbb{R}$ векторів $1,e_2,e_3$, які задовольняють
характеристичне рівняння (\ref{dopolnenije----2}). Тоді в алгебрі
$\mathbb{A}_{n-m+1}^1=1\oplus_s N$ (де нільпотентна підалгебра $N$
та ж сама що й в алгебрі $\mathbb{A}_n^m$) для кожного
$u\in\{1,2,\ldots,m\}$ існує трійка векторів
\begin{equation}\label{har-hv}
\begin{array}{l}
\widetilde{e}_1(u)=1,\vspace*{2mm} \\
\widetilde{e}_2(u):=a_u+I_u\,{\rm Rad}\,e_2\,,\vspace*{2mm} \\
\widetilde{e}_3(u):=b_u+I_u\,{\rm Rad}\,e_3\\
\end{array}
 \end{equation}
така, що характеристична система, породжена рівнянням
$\mathcal{X}(1,\widetilde{e}_2(u),\widetilde{e}_3(u))=0$ є редукцією
характеристичної системи, породженої рівнянням
(\ref{dopolnenije----2}).
\end{lemma}

\begin{proof}
Наслідком рівностей (\ref{har-hv}) є рівності
\begin{equation}\label{har-hv++}
\begin{array}{l}
\widetilde{e}_2^\beta(u)=a_u^\beta+I_u\sum\limits_{k=1}^\beta
\mathcal{C}_\beta^k\,a_u^{\beta-k} \left({\rm Rad}\,e_2\right)^k,\vspace*{2mm} \\
\widetilde{e}_3^\gamma(u)=b_u^\gamma+I_u\sum\limits_{k=1}^\gamma
\mathcal{C}_\gamma^k\,b_u^{\gamma-k} \left({\rm Rad}\,e_3\right)^k.
\end{array}
 \end{equation}
Враховуючи формули (\ref{har-hv++}), характеристичний многочлен
$\mathcal{X}(1,\widetilde{e}_2(u),\widetilde{e}_3(u))=0$ набуває
вигляду
$$
\sum\limits_{\alpha+\beta+\gamma=N}C_{\alpha,\beta,\gamma}\,\widetilde{e}_2^\beta(u)\,
\widetilde{e}_3^\gamma(u)=\sum\limits_{\alpha+\beta+\gamma=N}C_{\alpha,\beta,\gamma}\,\biggr(a_u^\beta\,b_u^\gamma+
$$
$$ +I_u\,b_u^\gamma\sum\limits_{k=1}^\beta
\mathcal{C}_\beta^k\,a_u^{\beta-k} \left({\rm
Rad}\,e_2\right)^k+I_u\,a_u^\beta\sum\limits_{k=1}^\gamma
\mathcal{C}_\gamma^k\,b_u^{\gamma-k} \left({\rm Rad}\,e_3\right)^k+
$$
\begin{equation}\label{dopolnenije----2+}
+I_u\sum\limits_{k=1}^\beta \mathcal{C}_\beta^k\,a_u^{\beta-k}
\left({\rm Rad}\,e_2\right)^k\sum\limits_{p=1}^\gamma
\mathcal{C}_\gamma^p\,b_u^{\gamma-p} \left({\rm
Rad}\,e_3\right)^p\biggr)=0.
\end{equation}

Далі покажемо, що характеристичні системи, породжені рівняннями
$\mathcal{X}(1,\widetilde{e}_2(u),\widetilde{e}_3(u))=0$ і
$\mathcal{X}(I_u,e_2I_u,e_3I_u)=0$ співпадають.

З цією метою зауважимо, що наслідком розкладів (\ref{e_1_e_2_e_3+})
є подання
$$e_2=a_1I_1+\cdots +a_mI_m+{\rm
Rad}\,e_2\,,\quad e_3=b_1I_1+\cdots +b_mI_m+{\rm Rad}\,e_3\,,
$$
з яких випливають співвідношення
\begin{equation}\label{har-hv+++}
e_2I_u=a_uI_u+I_u\,{\rm Rad}\,e_2\,,\quad e_3I_u=b_uI_u+I_u\,{\rm
Rad}\,e_3\,.
 \end{equation}

Тепер з (\ref{har-hv+++}) випливають рівності
\begin{equation}\label{har-hv++-}
\begin{array}{l}
e_2^\beta I_u=a_u^\beta I_u+I_u\sum\limits_{k=1}^\beta
\mathcal{C}_\beta^k\,a_u^{\beta-k} \left({\rm Rad}\,e_2\right)^k,\vspace*{2mm} \\
e_3^\gamma I_u=b_u^\gamma I_u+I_u\sum\limits_{k=1}^\gamma
\mathcal{C}_\gamma^k\,b_u^{\gamma-k} \left({\rm Rad}\,e_3\right)^k.
\end{array}
 \end{equation}
Беручи до уваги формули (\ref{har-hv++-}), характеристичне рівняння
$\mathcal{X}(I_u,e_2I_u,e_3I_u)=0$ набуває вигляду
$$
I_u\sum\limits_{\alpha+\beta+\gamma=N}C_{\alpha,\beta,\gamma}\,e_2^\beta\,
e_3^\gamma=\sum\limits_{\alpha+\beta+\gamma=N}C_{\alpha,\beta,\gamma}\,\biggr(a_u^\beta\,b_u^\gamma\,I_u+
$$
$$
+I_u\,b_u^\gamma\sum\limits_{k=1}^\beta
\mathcal{C}_\beta^k\,a_u^{\beta-k} \left({\rm
Rad}\,e_2\right)^k+I_u\,a_u^\beta\sum\limits_{k=1}^\gamma
\mathcal{C}_\gamma^k\,b_u^{\gamma-k} \left({\rm Rad}\,e_3\right)^k+
$$
\begin{equation}\label{dopolnenije----2+-}
+I_u\sum\limits_{k=1}^\beta \mathcal{C}_\beta^k\,a_u^{\beta-k}
\left({\rm Rad}\,e_2\right)^k\sum\limits_{p=1}^\gamma
\mathcal{C}_\gamma^p\,b_u^{\gamma-p} \left({\rm
Rad}\,e_3\right)^p\biggr)=0.
\end{equation}

З рівностей (\ref{dopolnenije----2+}), (\ref{dopolnenije----2+-})
очевидним чином випливає, що характеристичні системи, породжені
рівняннями $\mathcal{X}(1,\widetilde{e}_2(u),\widetilde{e}_3(u))=0$
і $\mathcal{X}(I_u,e_2I_u,e_3I_u)=0$ співпадають. Тепер доведення
леми випливає з леми \ref{lem-1-umb}.

\end{proof}

\begin{remark}\label{rem-1}
Відмітимо, що алгебра $\mathbb{A}_{n-m+1}^1=1\oplus_s N$ з базисом
$\{1,I_{m+1},\ldots,I_n\}$ є підалгеброю алгебри
$\mathbb{A}_n^m=S\oplus_s N$. Дійсно, будь-який елемент $a$ алгебри
$\mathbb{A}_n^m=S\oplus_s N$ вигляду
$$a=a_0I_1+a_0I_2+\cdots
+a_0I_m+a_{m+1}I_{m+1}+\cdots+a_nI_n=$$
$$=a_0(I_1+\cdots+I_m)+a_{m+1}I_{m+1}+\cdots+a_nI_n=a_0+a_{m+1}I_{m+1}+\cdots+a_nI_n$$ є представленням довільного
елементу алгебри $\mathbb{A}_{n-m+1}^1=1\oplus_s N$.
\end{remark}

З твердження \ref{Prop1} і леми  \ref{lem-2-umb} випливає така

\begin{theorem}\label{Th1}
Нехай в алгебрі $\mathbb{A}_n^m=S\oplus_s N$ існує трійка лінійно
незалежних над $\mathbb{R}$ векторів $1,e_2,e_3$, які задовольняють
характеристичне рівняння (\ref{dopolnenije----2}). Тоді в алгебрі
$\mathbb{A}_{n-m+1}^1=1\oplus_s N$ (де нільпотентна підалгебра $N$
та ж сама що й в алгебрі $\mathbb{A}_n^m$) для кожного
$u\in\{1,2,\ldots,m\}$ трійка векторів (\ref{har-hv})  задовольняє
характеристичне рівняння
$\mathcal{X}(1,\widetilde{e}_2(u),\widetilde{e}_3(u))=0$.
\end{theorem}

\begin{example}\label{Ex-1}
Розглянемо над полем $\mathbb{C}$  алгебру $\mathbb{A}_3^2$ з
таблицею множення (див., наприклад, \cite[c. 32]{Plaksa},
\cite{Pl-Pukh})
\begin{equation}\label{ex1}
\begin{tabular}{c||c|c||c|}
$\cdot$ & $I_1$ & $I_2$ & $I_3$ \\
\hline \hline
$I_1$ & $I_1$ & $0$ & $0$ \\
\hline
$I_2$ & $0$ & $I_2$ & $I_3$\\
\hline\hline
$I_3$ & $0$ & $I_3$ & $0$ \\
\hline
\end{tabular}\,\,.
\end{equation}
Очевидно, що напівпростою підалгеброю $S$ є підалгебра, породжена
ідемпотентами $I_1,I_2$, а нільпотентною підалгеброю $N$ є
підалгебра $\{\alpha I_3:\alpha\in\mathbb{C}\}$. Тоді алгебра
$\mathbb{A}_{2}^1:=1\oplus_s N$ співпадає з відомою бігармонічною
алгеброю $\mathbb{B}$ (див., наприклад, \cite{Gr-Pla-UMJ-2009}) і
має таку таблицю множення:
\begin{equation}\label{ex11}
\begin{tabular}{c||c|c|}
$\cdot$ & $1$ & $I_3$ \\
\hline \hline
$1$ & $1$ & $I_3$ \\
\hline
$I_3$ & $I_3$ & $0$ \\
\hline
\end{tabular}\,\,.
\end{equation}

Нехай в алгебрі $\mathbb{A}_3^2$ задане характеристичне рівняння
 (\ref{garmonichnyj_bazys-ogljad}).
Як відомо (див. теорему 1.8 в \cite{Plaksa}), умова
\textit{гармонічності} (\ref{garmonichnyj_bazys-ogljad}) векторів
$e_1=1$, $e_2=a_1I_1+a_2I_2+a_3I_3$, $e_3=b_1I_1+b_2I_2+b_3I_3$
алгебри $\mathbb{A}_3^2$ рівносильна системі рівнянь (\ref{s1}).

Оскільки для алгебри $\mathbb{A}_3^2$\, $m=2$, то в алгебрі
$\mathbb{B}$ ми будуємо дві трійки векторів виду (\ref{har-hv}):
\begin{equation}\label{har-hv=}
\widetilde{e}_1(1)=1,\,\widetilde{e}_2(1)=a_1+I_1(a_3I_3)=a_1\,,\,
\widetilde{e}_3(1)=b_1+I_1(b_3I_3)=b_1\\
\end{equation}
та
\begin{equation}\label{har-hv==}
\begin{array}{l}
\widetilde{e}_1(2)=1,\vspace*{2mm} \\
\widetilde{e}_2(2)=a_2+I_2(a_3I_3)=a_2+a_3I_3\,,\vspace*{2mm} \\
\widetilde{e}_3(2)=b_2+I_2(b_3I_3)=b_2+b_3I_3\,.\\
\end{array}
 \end{equation}
За теоремою \ref{Th1} трійки (\ref{har-hv=}) та (\ref{har-hv==})
гармонічні в алгебрі $\mathbb{B}$ (тобто задовольняють умову
(\ref{garmonichnyj_bazys-ogljad})). Справді, гармонічність трійки
(\ref{har-hv=}) рівносильна першому рівнянню системи (\ref{s1}), а
гармонічність трійки (\ref{har-hv==}) рівносильна системі
(\ref{s2}).
\end{example}

\begin{example}\label{Ex-2}
Розглянемо над полем $\mathbb{C}$  алгебру $\mathbb{A}_5^3$ з такою
таблицею множення
\begin{equation}\label{ex2}
\begin{tabular}{c||c|c|c||c|c|}
$\cdot$ & $I_1$ & $I_2$ & $I_3$& $I_4$ & $I_5$  \\
\hline \hline
$I_1$ & $I_1$ & $0$ & $0$ & $0$& $I_5$\\
\hline
$I_2$ & $0$ & $I_2$ & $0$& $0$& $0$\\
\hline
$I_3$ & $0$ & $0$ & $I_3$ & $I_4$& $0$\\
\hline\hline
$I_4$ & $0$ & $0$ & $I_4$ & $0$& $0$\\
\hline
$I_5$ & $I_5$ & $0$ & $0$& $0$& $0$ \\
\hline\hline
\end{tabular}\,\,.
\end{equation}
Відмітимо, що напівпростою підалгеброю $S$ є підалгебра, породжена
ідемпотентами $I_1,I_2,I_3$, а нільпотентною підалгеброю $N$ є
підалгебра з базисом $\{I_4,I_5\}$. Тоді алгебра
$\mathbb{A}_{3}^1:=1\oplus_s N$ співпадає з відомою алгеброю
$\mathbb{A}_4$ (див., наприклад, \cite[c. 26]{Plaksa}) і має таку
таблицю множення:
\begin{equation}\label{ex21}
\begin{tabular}{c||c|c|c|}
$\cdot$ & $1$ & $I_4$ & $I_5$ \\
\hline \hline
$1$ & $1$ & $I_4$  & $I_5$\\
\hline
$I_4$ & $I_4$ & $0$ & $0$ \\
\hline
$I_5$ & $I_5$ & $0$  & $0$\\
\hline
\end{tabular}\,\,.
\end{equation}

Нехай в алгебрі $\mathbb{A}_5^3$ задане характеристичне рівняння
 (\ref{garmonichnyj_bazys-ogljad}).
Умова гармонічності (\ref{garmonichnyj_bazys-ogljad}) векторів
вигляду (\ref{e_1_e_2_e_3+}) алгебри $\mathbb{A}_5^3$ рівносильна
наступній системі рівнянь
\begin{equation}\label{har-hv==-}
\begin{array}{l}
1+a_u^2+b_u^2=0,\quad u=1,2,3,\vspace*{2mm} \\
a_3a_4+b_3b_4=0,\vspace*{2mm} \\
a_1a_5+b_1b_5=0.\\
\end{array}
 \end{equation}

 Оскільки для алгебри $\mathbb{A}_5^3$\, $m=3$, то
в алгебрі $\mathbb{A}_4$ ми будуємо три трійки векторів виду
(\ref{har-hv}):
\begin{equation}\label{har-hv==+}
\begin{array}{l}
\widetilde{e}_1(1)=1,\vspace*{2mm} \\
\widetilde{e}_2(1)=a_1+I_1(a_4I_4+a_5I_5)=a_1+a_5I_5\,,\vspace*{2mm} \\
\widetilde{e}_3(1)=b_1+I_1(b_4I_4+b_5I_5)=b_1+b_5I_5\,,\\
\end{array}
 \end{equation}
\begin{equation}\label{har-hv==++}
\begin{array}{l}
\widetilde{e}_1(2)=1,\vspace*{2mm} \\
\widetilde{e}_2(2)=a_2+I_2(a_4I_4+a_5I_5)=a_2\,,\vspace*{2mm} \\
\widetilde{e}_3(2)=b_2+I_2(b_4I_4+b_5I_5)=b_2\,,\\
\end{array}
 \end{equation}
та
\begin{equation}\label{har-hv==+++}
\begin{array}{l}
\widetilde{e}_1(3)=1,\vspace*{2mm} \\
\widetilde{e}_2(3)=a_3+I_3(a_4I_4+a_5I_5)=a_3+a_4I_4\,,\vspace*{2mm} \\
\widetilde{e}_3(3)=b_3+I_3(b_4I_4+b_5I_5)=b_3+b_4I_4\,.\\
\end{array}
\end{equation}

За теоремою \ref{Th1} трійки (\ref{har-hv==+}), (\ref{har-hv==++})
та (\ref{har-hv==+++}) гармонічні в алгебрі $\mathbb{A}_4$ (тобто
задовольняють умову (\ref{garmonichnyj_bazys-ogljad})). Справді,
гармонічність трійки (\ref{har-hv==+}) рівносильна системі з першого
і п'ятого рівняння системи (\ref{har-hv==}); гармонічність трійки
(\ref{har-hv==++}) рівносильна другому рівнянню системи
(\ref{har-hv==}), а гармонічність трійки (\ref{har-hv==+++})
рівносильна системі з третього і четвертого рівняння системи
(\ref{har-hv==}).

\end{example}

\subsection{Лінійна незалежність векторів $1,\widetilde{e}_2(u),\widetilde{e}_3(u)$}

З наведених прикладів видно, що вектори
$1,\widetilde{e}_2(u),\widetilde{e}_3(u)$ при деяких
$u\in\{1,2,\ldots,m\}$ можуть бути лінійно залежними над полем
$\mathbb{R}$. Так, трійки (\ref{har-hv=}), (\ref{har-hv==++}) завжди
лінійно залежні над полем $\mathbb{R}$.

Встановимо необхідні і достатні умови лінійної незалежності над
полем $\mathbb{R}$ векторів
$1,\widetilde{e}_2(u),\widetilde{e}_3(u)$ алгебри
$\mathbb{A}_{n-m+1}^1=1\oplus_s N$.

\begin{lemma}\label{lem-3-umb}
Нехай вектори (\ref{e_1_e_2_e_3+}) алгебри $\mathbb{A}_n^m=S\oplus_s
N$ лінійно незалежні над полем $\mathbb{R}$ і нехай
$u\in\{1,2,\ldots,m\}$ фіксоване. Тоді

1. якщо вектори $I_u\,{\rm Rad}\,e_2$, $I_u\,{\rm
Rad}\,e_3\in\mathbb{A}_n^m$ лінійно незалежні над полем
$\mathbb{R}$, то вектори $1,\widetilde{e}_2(u),\widetilde{e}_3(u)$
алгебри $\mathbb{A}_{n-m+1}^1=1\oplus_s N$ також лінійно незалежні
над полем $\mathbb{R}$;

2. якщо ж вектори $I_u\,{\rm Rad}\,e_2\,$, $I_u\,{\rm
Rad}\,e_3\in\mathbb{A}_n^m$ лінійно залежні над полем $\mathbb{R}$,
то вектори $1,\widetilde{e}_2(u),\widetilde{e}_3(u)$ алгебри
$\mathbb{A}_{n-m+1}^1=1\oplus_s N$ лінійно незалежні над полем
$\mathbb{R}$ тоді і тільки тоді, коли існує $r\in\{m+1,\ldots,n\}$
таке, що $I_uI_r=I_r$ і виконується хоча б одне співвідношення
\begin{equation}\label{har-riv-+}
{\rm Im}\,a_u \,{\rm Re}\,b_r \neq{\rm Im}\,b_u\, {\rm Re}\,a_r
\quad \text{або} \quad {\rm Im}\,a_u \,{\rm Im}\,b_r \neq{\rm
Im}\,b_u\, {\rm Im}\,a_r \,.
\end{equation}
\end{lemma}

\begin{proof}
Доведемо перше твердження леми. За умовою рівність
\begin{equation}\label{har-riv+}
\beta_2\,I_u\,{\rm Rad}\,e_2+\beta_3\,I_u\,{\rm Rad}\,e_3=0,\quad
\beta_2,\beta_3\in\mathbb{R}
\end{equation}
виконується тоді і тільки тоді, коли $\beta_2=\beta_3=0$.

Розглянемо лінійну комбінацію
$$
\alpha_1+\alpha_2\,\widetilde{e}_2(u)+\alpha_3\,\widetilde{e}_3(u)=
(\alpha_1+\alpha_2\,a_u+\alpha_3\,b_u)+
$$
\begin{equation}\label{har-riv+1}
+\left(\alpha_2\,I_u\,{\rm Rad}\,e_2+\alpha_3\,I_u\,{\rm
Rad}\,e_3\right)=0, \quad \alpha_1,\alpha_2,\alpha_3\in\mathbb{R}.
\end{equation}
Оскільки вираз у другій дужці в рівності (\ref{har-riv+1}) приймає
значення в радикалі $\mathcal{R}$ алгебри, а перша дужка
комплекснозначна, то умова (\ref{har-riv+1}) рівносильна системі
рівнянь
\begin{equation}\label{har-riv+2}
\begin{array}{l}
\alpha_1+\alpha_2\,a_u+\alpha_3\,b_u=0,\vspace*{2mm} \\
\alpha_2\,I_u\,{\rm
Rad}\,e_2+\alpha_3\,I_u\,{\rm Rad}\,e_3=0.\\
\end{array}
\end{equation}
З другого рівняння системи (\ref{har-riv+2}) і умови
(\ref{har-riv+}) випливає, що $\alpha_2=\alpha_3=0$. А тоді з
першого рівняння системи (\ref{har-riv+2}) отримуємо $\alpha_1=0$.
Отже, вектори $1,\widetilde{e}_2(u),\widetilde{e}_3(u)$ лінійно
незалежні над $\mathbb{R}$.

Доведемо друге твердження леми. Розглянемо рівність
$$
\beta_1+\beta_2\,e_2+\beta_3\,e_3=\sum\limits_{s=1}^m
I_s(\beta_1+\beta_2\,a_s+\beta_3\,b_s)+ \sum\limits_{k=m+1}^n
I_k(\beta_2\,a_k+\beta_3\,b_k)=0,
$$
яка рівносильна системі рівнянь
\begin{equation}\label{har-riv+3}
\begin{array}{rcl}
\beta_1+\beta_2\,{\rm Re}\,a_s+\beta_3\,{\rm Re}\,b_s&=&0,\vspace*{2mm} \\
\beta_2\,{\rm Im}\,a_s+\beta_3\,{\rm Im}\,b_s&=&0,\qquad s=1,2,\ldots,m,\vspace*{2mm} \\
\beta_2\,{\rm Re}\,a_k+\beta_3\,{\rm Re}\,b_k&=&0,\vspace*{2mm} \\
\beta_2\,{\rm Im}\,a_k+\beta_3\,{\rm Im}\,b_k&=&0,\qquad k=m+1,\ldots,n.\\
\end{array}
\end{equation}
Лінійна незалежність над $\mathbb{R}$ векторів $1,e_2,e_3$ означає,
що серед усіх рівнянь системи (\ref{har-riv+3}), окрім першого,
існує хоча б два рівняння, які між собою не пропорційні.

Тепер запишемо умову лінійної незалежності над $\mathbb{R}$ векторів
$1,\widetilde{e}_2(u),\widetilde{e}_3(u)$. Для цього систему
(\ref{har-riv+2}) запишемо в розгорнутому вигляді
\begin{equation}\label{har-riv+4}
\begin{array}{rcl}
\alpha_1+\alpha_2\,{\rm Re}\,a_u+\alpha_3\,{\rm Re}\,b_u&=&0,\vspace*{2mm} \\
\alpha_2\,{\rm Im}\,a_u+\alpha_3\,{\rm Im}\,b_u&=&0,\vspace*{2mm} \\
\alpha_2\,{\rm Re}\,a_r+\alpha_3\,{\rm Re}\,b_r&=&0,\vspace*{2mm} \\
\alpha_2\,{\rm Im}\,a_r+\alpha_3\,{\rm Im}\,b_r&=&0\vspace*{2mm}\\
 \forall\, r\in\{m+1,\ldots,n\} &:& I_uI_r=I_r\,.\\
\end{array}
\end{equation}

За умовою пункту \textit{2} леми вектори $I_u\,{\rm Rad}\,e_2\,$,
$I_u\,{\rm Rad}\,e_3$ лінійно залежні над $\mathbb{R}$. Це означає,
що в системі (\ref{har-riv+4}) всі рівності, окрім перших двох,
пропорційні між собою. Очевидно, що для лінійної незалежності над
$\mathbb{R}$ векторів $1,\widetilde{e}_2(u),\widetilde{e}_3(u)$
необхідно і достатньо, щоб друге рівняння системи (\ref{har-riv+4})
було не пропорційне хоча б з одним іншим рівнянням (крім першого)
системи (\ref{har-riv+4}). А це рівносильно умовам
(\ref{har-riv-+}).

\end{proof}

\section{Моногенні функції, визначені в різних комутативних алгебрах}

В алгебрі $\mathbb{A}_n^m=S\oplus_s N$ будемо розглядати моногенні
функції $\Phi$, визначені в деякій області $\Pi_\zeta\subset E_3$
виду (\ref{Pi}). Геометрично область $\Pi\subset\mathbb{R}^3$, яка
конгуентна області $\Pi_\zeta\subset E_3$, є перетином $m$
нескінченних циліндрів, кожен з яких паралельний деякій з $m$ прямих
$L_u$\,, $u=1,2,\ldots,m$ вигляду (\ref{L-u}). Тобто,
$\Pi=\cap_{u=1}^m\Pi(u),$ де $\mathbb{R}^3\supset\Pi(u)$
--- нескінченний циліндр, паралельний прямій $L_u$\,. І те ж саме маємо для
конгруентних областей в $E_3$:
\begin{equation}\label{p5-0}
\Pi_\zeta=\bigcap\limits_{u=1}^m\Pi_\zeta(u).
\end{equation}
Аналітично циліндр $\Pi_\zeta(u)$ визначається рівністю
$$
\Pi_\zeta(u)=\{\zeta_u:=I_u\,\zeta \,:\, \zeta\in\Pi_\zeta\}.
$$

Тепер розглянемо моногенну в області  $\Pi_\zeta$ функцію
$\Phi:\Pi_\zeta\rightarrow\mathbb{A}_n^m$. Введемо позначення
\begin{equation}\label{p5}
\Phi_u(\zeta):=I_u\,\Phi(\zeta),\quad u=1,2,\ldots,m.
\end{equation}
Тоді очевидною є рівність
\begin{equation}\label{p5-1}
\Phi=(I_1+\cdots +I_m)\Phi=\sum_{u=1}^m\Phi_u\,.
\end{equation}
Крім того, з рівності (\ref{Teor--1}) і таблиці множення алгебри
$\mathbb{A}_n^m$ випливає, що при кожному $u\in\{1,2,\ldots,m\}$
функція $\Phi_u$ моногенна у всьому нескінченному циліндрі
$\Pi_\zeta(u)$.

Таким чином, кожна моногенна в області (\ref{p5-0}) функція
$\Phi:\Pi_\zeta\rightarrow\mathbb{A}_n^m$ подається у вигляді суми
(\ref{p5-1}), де функція $\Phi_u$ моногенна у всьому циліндрі
$\Pi_\zeta(u)$.

Тепер перейдемо до розгляду моногенних функцій $\widetilde{\Phi}$ в
алгебрі $\mathbb{A}_{n-m+1}^1=1\oplus_s N$. Оскільки, відповідно до
зауваження \ref{rem-1}, алгебра $\mathbb{A}_{n-m+1}^1$ є підалгеброю
алгебри $\mathbb{A}_n^m$, то в алгебрі $\mathbb{A}_{n-m+1}^1$ усі
циліндри $\Pi_\zeta(u)$ з рівності (\ref{p5-0}) співпадають між
собою. Тобто, в алгебрах вигляду $\mathbb{A}_{n-m+1}^1$ кожна
моногенна функція буде моногенною в деякому одному нескінченному
циліндрі.

В наступній теоремі встановлюється зв'язок між моногенними функціями
в алгебрах $\mathbb{A}_n^m=S\oplus_s N$ та
$\mathbb{A}_{n-m+1}^1=1\oplus_s N$. Для формулювання результату
введемо деякі позначення.

На вектори вигляду (\ref{har-hv}) алгебри $\mathbb{A}_{n-m+1}^1$
натягнемо лінійний простір
$\widetilde{E}_3(u):=\{\widetilde{\zeta}(u)=x+y\widetilde{e}_2(u)+z\widetilde{e}_3(u):\,\,
x,y,z\in\mathbb{R}\}$. Трійка векторів (\ref{har-hv}) визначає одну
пряму $\widetilde{L}(u)$ виду (\ref{L-u}), яка відповідає множині
необоротних елементів $\widetilde{\zeta}(u)$ простору
$\widetilde{E}_3(u)$. Нехай $\widetilde{\Pi}_{\widetilde{\zeta}(u)}$
--- деякий нескінченний циліндр в
$\widetilde{E}_3(u)$, паралельний прямій $\widetilde{L}(u)$.

\begin{theorem}\label{Teor--2}
Нехай в алгебрі $\mathbb{A}_n^m=S\oplus_s N$ існує трійка лінійно
незалежних над $\mathbb{R}$ векторів $1,e_2,e_3$, які задовольняють
характеристичне рівняння (\ref{dopolnenije----2}) і нехай
 $f_u(E_3)=\mathbb{C}$ при всіх
 $u=1,2,\ldots, m$. Крім того, нехай функція $\Phi:\Pi_\zeta\rightarrow\mathbb{A}_n^m$ змінної
$\zeta=x+ye_2+ze_3$ моногенна в області $\Pi_\zeta\subset E_3$ виду
(\ref{p5-0}). Тоді в алгебрі $\mathbb{A}_{n-m+1}^1=1\oplus_s N$ (де
нільпотентна підалгебра $N$ та ж сама що й в алгебрі
$\mathbb{A}_n^m$) для кожного $u\in\{1,2,\ldots,m\}$ існує трійка
векторів (\ref{har-hv}), яка задовольняє характеристичне рівняння
$\mathcal{X}(1,\widetilde{e}_2(u),\widetilde{e}_3(u))=0$ і існує
функція
$\widetilde{\Phi}_u:\widetilde{\Pi}_{\widetilde{\zeta}(u)}\rightarrow\mathbb{A}_{n-m+1}^1$
змінної $\widetilde{\zeta}(u)$, яка моногенна в циліндрі
$$
\widetilde{\Pi}_{\widetilde{\zeta}(u)}=\left\{\widetilde{\zeta}(u)\in\widetilde{E}_3(u)\,:\,
f_u\big(\widetilde{\zeta}(u)\big)=f_u(\zeta)\,,\,\zeta\in
\Pi_\zeta(u)\right\}
$$
 така, що
\begin{equation}\label{p5-2}
\Phi_u(\zeta)=I_u\,\widetilde{\Phi}_u\big(\widetilde{\zeta}(u)\big).
\end{equation}
\end{theorem}

\begin{proof} Існування трійки (\ref{har-hv}) з властивістю
$\mathcal{X}(1,\widetilde{e}_2(u),\widetilde{e}_3(u))=0$ доведено в
теоремі \ref{Th1}. Нехай надалі $u\in\{1,2,\ldots,m\}$ фіксоване.
Доведемо існування і моногенність в області
$\widetilde{\Pi}_{\widetilde{\zeta}(u)}$ функції
$\widetilde{\Phi}_u$, яка задовольняє рівність (\ref{p5-2}). З цією
метою спочатку доведемо рівність
\begin{equation}\label{p5-3}
I_u\,\zeta^{-1}=I_u\,\widetilde{\zeta}^{-1}(u)
\end{equation}
$$\forall\,\zeta=x+ye_2+ze_3 \,\,\,
\forall\,\widetilde{\zeta}(u)=x+y\widetilde{e}_2(u)+z\widetilde{e}_3(u),\,\,x\in\mathbb{C},y,z\in\mathbb{R}.
$$

 З рівностей (\ref{har-hv++-}),  (\ref{har-hv})  випливають
співвідношення
$$
I_u\,e_2=I_u\,\widetilde{e}_2(u)\,,\quad
I_u\,e_3=I_u\,\widetilde{e}_3(u)\,,
$$
з яких, в свою чергу, випливає рівність
\begin{equation}\label{p5-4}
I_u\,\zeta=I_u\,\widetilde{\zeta}(u).
\end{equation}

Розглянемо різницю $I_u\,\zeta^{-1}-I_u\,\widetilde{\zeta}^{-1}(u)$.
За формулою Гільберта (див., наприклад, теорему 4.8.2 в
\cite{Hil_Filips}), маємо
$$
I_u\,\zeta^{-1}-I_u\,\widetilde{\zeta}^{-1}(u)=\big(I_u\,\zeta-I_u\,
\widetilde{\zeta}(u)\big)\Big(\zeta\,\widetilde{\zeta}(u)\Big)^{-1}=0,
$$
внаслідок рівності (\ref{p5-4}). Отже, рівність (\ref{p5-3})
доведено. Тепер з (\ref{p5-3}) маємо співвідношення
\begin{equation}\label{p5-5}
I_u(t-\zeta)^{-1}=I_u\big(t-\widetilde{\zeta}(u)\big)^{-1}
\end{equation}
$$ \forall\, t\in\mathbb{C}:\, t\neq\xi_u=f_u(\zeta)\quad
\forall\,\zeta\in\Pi_\zeta(u)\quad\forall\,\widetilde{\zeta}(u)\in\widetilde{\Pi}_{\widetilde{\zeta}(u)}\,.
$$

З таблиці множення алгебри $\mathbb{A}_n^m$ і формули
(\ref{Teor--1}) для моногенної в області $\Pi_\zeta(u)$ функції
$\Phi_u(\zeta)$ маємо представлення
 \begin{equation}\label{p5-6}
\Phi_u(\zeta)=I_u\,\frac{1}{2\pi i}\int\limits_{\Gamma_u}
\Bigr(F_u(t)+
\sum\limits_{s=m+1}^nI_s\,G_s(t)\Bigr)(t-\zeta)^{-1}\,dt,
 \end{equation}
де функції $F_u\,,G_s$ визначені в теоремі $\textbf{A}$.

Враховуючи співвідношення (\ref{p5-5}), представлення (\ref{p5-6})
перепишемо у вигляді
\begin{equation}\label{p5-7}
\Phi_u(\zeta)=I_u\,\frac{1}{2\pi i}\int\limits_{\Gamma_u}
\Bigr(F_u(t)+
\sum\limits_{s=m+1}^nI_s\,G_s(t)\Bigr)\big(t-\widetilde{\zeta}(u)\big)^{-1}\,dt.
\end{equation}

Оскільки в алгебрі $\mathbb{A}_{n-m+1}^1$ міститься єдиний
максимальний ідеал $\mathcal{I}$, який співпадає з радикалом
(\ref{radyk}) цієї алгебри $\mathcal{R}$, то на цій алгебрі
визначений єдиний лінійний неперервний мультиплікативний функціонал
$f:\mathbb{A}_{n-m+1}^1\rightarrow\mathbb{C}$ ядром якого є радикал
$\mathcal{R}$. А це означає, що
$f\big(\widetilde{\zeta}(u)\big)=x+a_u\,y+b_u\,z$ для кожного
$\widetilde{\zeta}(u)\in\widetilde{E}_3(u)$. Беручи до уваги
рівність $f_u(\zeta)=x+a_u\,y+b_u\,z$ для довільного $\zeta\in E_3$,
маємо рівність
\begin{equation}\label{p5-8}
f\big(\widetilde{\zeta}(u)\big)=f_u(\zeta).
\end{equation}
З рівності (\ref{p5-8}) і умови теореми $f_u(E_3)=\mathbb{C}$
 отримуємо співвідношення
 $f\big(\widetilde{\zeta}(u)\big)=\mathbb{C}$ для довільного
$\widetilde{\zeta}(u)\in\widetilde{E}_3(u)$.

Щойно ми показали, що виконуються умови теореми $\textbf{A}$ для
моногенних функцій в алгебрі $\mathbb{A}_{n-m+1}^1$. Тоді в алгебрі
$\mathbb{A}_{n-m+1}^1$ формула (\ref{Teor--1}) для моногенної в
області $\widetilde{\Pi}_{\widetilde{\zeta}(u)}$ функції
$\widetilde{\Phi}_u\big(\widetilde{\zeta}(u)\big)$ має вигляд
 \begin{equation}\label{p5-9}
\widetilde{\Phi}_u\big(\widetilde{\zeta}(u)\big)=\frac{1}{2\pi
i}\int\limits_{\gamma} \Bigr(\widetilde{F}(t)+
\sum\limits_{s=m+1}^nI_s\,\widetilde{G}_s(t)\Bigr)\big(t-\widetilde{\zeta}(u)\big)^{-1}\,dt.
 \end{equation}

Потрібна нам формула (\ref{p5-2}) буде прямим наслідком рівностей
(\ref{p5-7}) та (\ref{p5-9}), якщо ми покажемо, що можна покласти
$\gamma\equiv\Gamma_u$, $F_u\equiv\widetilde{F}$ і
$G_s\equiv\widetilde{G}_s$ для таких $s$, що $I_u\,I_s=I_s$\,.
Покажемо це.

З рівності (\ref{p5-8}) випливає, що циліндри $\Pi_\zeta(u)\subset
E_3$ та $\widetilde{\Pi}_{\widetilde{\zeta}(u)}\subset
\widetilde{E}_3(u)$ відповідними функціоналами $f_u$ та $f$
відображаються в одну і ту ж область $D$ комплексної площини
$\mathbb{C}$. А це означає, що функції $F_u$ і $\widetilde{F}$, а
також функції $G_s$ і $\widetilde{G}_s$ голоморфні в одній і тій
самій області $D$. Отже, ми можемо покласти $F_u\equiv\widetilde{F}$
і $G_s\equiv\widetilde{G}_s$ в $D$.

Оскільки криві інтегрування $\gamma$ і $\Gamma_u$ лежать в області
$D$, то ми можемо взяти $\gamma\equiv\Gamma_u$\,. Більше того,
оскільки за теоремою $\textbf{А}$ крива $\Gamma_u$ в рівності
(\ref{p5-6}) охоплює точку $f_u(\zeta)=x+a_u\,y+b_u\,z$, то
внаслідок рівності (\ref{p5-8}) крива $\gamma\equiv\Gamma_u$ охоплює
спектр точки $\widetilde{\zeta}(u)$
--- точку $f\big(\widetilde{\zeta}(u)\big)=x+a_u\,y+b_u\,z$. А це
нам і потрібно. Теорему доведено.
\end{proof}

\begin{remark}
З рівностей (\ref{p5-1}), (\ref{p5-2}) випливає представлення
\begin{equation}\label{p5-10}
\Phi(\zeta)=I_1\,\widetilde{\Phi}_1\big(\widetilde{\zeta}(1)\big)+\cdots
+I_m\,\widetilde{\Phi}_m\big(\widetilde{\zeta}(m)\big).
\end{equation}
\end{remark}

\begin{remark}
Теорема \ref{Th1} означає, що функції $\Phi$ і
$I_u\widetilde{\Phi}_u$ при всіх $u=1,2,\ldots,m$ задовольняють одне
й те ж саме диференціальне рівняння виду (\ref{dopolnenije----1}).
\end{remark}

\begin{remark}
Теорема \ref{Teor--2} стверджує, що для побудови розв'язків
диференціального рівняння (\ref{dopolnenije----1}) у вигляді
компонент моногенних функцій зі значеннями в комутативних алгебрах,
достатньо обмежитись вивченням моногенних функцій в алгебрах з
базисом $\{1,\eta_1,\eta_2,\ldots,\eta_n\}$, де
$\eta_1,\eta_2,\ldots,\eta_n$ --- нільпотенти. Тобто кількість таких
$n$-вимірних комутативних асоціативних алгебр з одиницею над полем
$\mathbb{C}$ в яких потрібно вивчати моногенні функції рівна
кількості $(n-1)$-вимірних комутативних асоціативних комплексних
нільпотентних алгебр.

Зокрема, серед двовимірних комутативних асоціативних алгебр з
одиницею над полем $\mathbb{C}$ (яких існує всього дві) достатньо
обмежитись вивченням моногенних функцій в бігармонічній алгебрі
$\mathbb{B}$. Серед тривимірних комутативних асоціативних алгебр з
одиницею над полем $\mathbb{C}$ (яких існує всього чотири) достатньо
обмежитись вивченням моногенних функцій в двох із них (це алгебри
$\mathbb{A}_3$ і $\mathbb{A}_4$ в термінах роботи \cite{Plaksa}). А
серед чотиривимірних комутативних асоціативних алгебр з одиницею над
полем $\mathbb{C}$ (яких існує всього 9, див. \cite{Mazzola-80})
достатньо обмежитись вивченням моногенних функцій в чотирьох із них
(це алгебри $\widetilde{A}_{3,1}$, $\widetilde{A}_{3,2}$,
$\widetilde{A}_{3,3}$, $\widetilde{A}_{3,4}$ з таблиці 9 роботи
\cite{Burde_de_Graaf}, див. також теорему 5.1 в роботі
\cite{Hegazim-Ab}). Серед усіх п'ятивимірних комутативних
асоціативних алгебр з одиницею над полем $\mathbb{C}$ (яких існує
всього 25, див. \cite{Mazzola-80}) достатньо обмежитись вивченням
моногенних функцій в дев'яти із них (таблиці множення усіх цих 9
нільпотентних чотиривимірних алгебр наведено в теоремі 6.1 з роботи
\cite{Hegazim-Ab}). І нарешті серед усіх шестивимірних комутативних
асоціативних алгебр з одиницею над полем $\mathbb{C}$ достатньо
обмежитись вивченням моногенних функцій в 25-ти із них (усі ці 25
нільпотентних п'ятивимірних алгебр наведено в таблиці 1 з роботи
\cite{Poonen}). Відомо також (див. \cite{Suprunenko}), що починаючи
з розмірності 6 множина усіх попарно неізоморфних нільпотентних
комутативних алгебр над $\mathbb{C}$  є нескінченною.
\end{remark}

\begin{remark}
Теорема \ref{Teor--2} залишається справедливою для випадку, коли ми
будемо розглядати функції $\Phi:\Pi_\zeta\rightarrow\mathbb{A}_n^m$
змінної $\zeta:=\sum\limits_{r=1}^kx_re_r$, $2\leq k\leq2n$, яка
моногенна в області $\Pi_\zeta\subset E_k$. При цьому замість
теореми $\textbf{A}$ необхідно використовувати теорему 1 з роботи
\cite{Shpakivskyi-Zb-2015-1}.
\end{remark}

Продемоструємо теорему \ref{Teor--2} на алгебрах, які розглядалися в
прикладах \ref{Ex-1} та  \ref{Ex-2}.

\begin{example}
Отже, розглядаємо алгебру $\mathbb{A}_3^2$ з таблицею множення
(\ref{ex1}). Для алгебри $\mathbb{A}_3^2$ алгеброю виду $1\oplus_s
N$ є бігармонічні алгебра $\mathbb{B}$ з таблицею множення
(\ref{ex11}).

Відповідно до представлення (\ref{Teor--1}), кожна моногенна функція
$\Phi$ зі значеннями в алгебрі $\mathbb{A}_3^2$ подається у вигляді
\begin{equation}\label{p5-11}
\Phi(\zeta)=F_1(\xi_1)I_1+F_2(\xi_2)I_2+\Big((a_3y+b_3z)F_2'(\xi_2)+G_3(\xi_2)\Big)I_3
\end{equation}
$$\forall\,\zeta\in\Pi_\zeta\,,\quad \xi_u=x+a_uy+b_uz,\,\,u=1,2,
$$
де $F_1$ --- деяка голоморфна функція в області $D_1$\,, а
$F_2\,,G_3$
--- деякі голоморфні функції в області $D_2$\,. Оскільки в
$\mathbb{A}_3^2$ \, $m=2$, то геометрично область $\Pi_\zeta$ є
перетином двох нескінченних циліндрів:\,
$\Pi_\zeta=\Pi_\zeta(1)\cap\Pi_\zeta(2)$.

Зауважимо, що представлення  (\ref{p5-11}) раніше було отримано в
роботі \cite{Pl-Pukh}. Крім того, функція (\ref{p5-11}) задовольняє
деяке диференціальне рівняння вигляду (\ref{dopolnenije----1}).

Подамо функцію (\ref{p5-11}) у вигляді (\ref{p5-1}):
\begin{equation}\label{p5-13}
\Phi(\zeta)=\Phi(\zeta)I_1+\Phi(\zeta)I_2=:\Phi_1(\zeta)+\Phi_2(\zeta),
\end{equation}
де $\Phi_1(\zeta)=F_1(\xi_1)I_1$ --- моногенна функція в циліндрі
$\Pi_\zeta(1)$, а функція
$$\Phi_2(\zeta)=F_2(\xi_2)I_2+\Big((a_3y+b_3z)F_2'(\xi_2)+G_3(\xi_2)\Big)I_3$$
моногенна в циліндрі $\Pi_\zeta(2)$.\vskip2mm

Перейдемо до розгляду моногенних функцій в алгебрі $\mathbb{B}$.  З
представлення (\ref{Teor--1}) випливає, що кожна моногенна функція
$\widetilde{\Phi}$ зі значеннями в алгебрі $\mathbb{B}$ подається у
вигляді
\begin{equation}\label{p5-12}
\widetilde{\Phi}(\widetilde{\zeta})=\widetilde{F}(\widetilde{\xi})+
\Big((a_3y+b_3z)\widetilde{F}'(\widetilde{\xi})+\widetilde{G}(\widetilde{\xi})\Big)I_3
\quad\forall\,\widetilde{\zeta}\in\widetilde{\Pi}_{\widetilde{\zeta}}\,,\,\,\widetilde{\xi}=f(\widetilde{\zeta}),
\end{equation}
де $\widetilde{F},\widetilde{G}$
--- деякі голоморфні функції в області $D$. Область $\widetilde{\Pi}_{\widetilde{\zeta}}$ є
нескінченним циліндром. Рівність (\ref{p5-12}) для спеціального
випадку встановлена в роботі \cite{Gr-Pla-UMJ-2009}.

Теорема \ref{Teor--2} стверджує наступне:

$\displaystyle 1)$ в алгебрі $\mathbb{B}$ існує трійка векторів
$1,\widetilde{e}_2(1)\,,\widetilde{e}_3(1)$, яка задовольняє те ж
саме характеристичне рівняння що й трійка
$1,e_2,e_3\in\mathbb{A}_3^2$. При цьому, будуть виконуватись
співвідношення $\xi_1\equiv\widetilde{\xi}$, $D_1\equiv D$ і, крім
того, існує моногенна в $\mathbb{B}$ функція $\widetilde{\Phi}$
така, що
\begin{equation}\label{p5-14}
I_1\,\widetilde{\Phi}_1\big(\widetilde{\zeta}(1)\big)=\Phi_1(\zeta)\,.
\end{equation}

$\displaystyle 2)$ в алгебрі $\mathbb{B}$ існує трійка векторів
$1,\widetilde{e}_2(2)\,,\widetilde{e}_3(2)$, яка задовольняє те ж
саме характеристичне рівняння що й трійка
$1,e_2,e_3\in\mathbb{A}_3^2$. При цьому, будуть виконуватись
співвідношення $\xi_2\equiv\widetilde{\xi}$, $D_2\equiv D$ і, крім
того, існує моногенна функція $\widetilde{\Phi}$ така, що
\begin{equation}\label{p5-15}
I_2\,\widetilde{\Phi}_2\big(\widetilde{\zeta}(2)\big)=\Phi_2(\zeta)\,.
\end{equation}

Потрібні трійки векторів $1,\widetilde{e}_2(1)\,,\widetilde{e}_3(1)$
та $1,\widetilde{e}_2(2)\,,\widetilde{e}_3(2)$ були знайдені у
прикладі \ref{Ex-1}. Розглянемо випадок 1). Дійсно, для трійки
(\ref{har-hv=}) маємо
$\widetilde{\zeta}(1)=x+a_1y+b_1z\equiv\xi_1\equiv\widetilde{\xi}$,
$D_1\equiv D$. Покладемо $\widetilde{F}\equiv F_1$,
$\widetilde{G}\equiv G_3$ в $D$. Тоді рівність (\ref{p5-12})
перепишеться у вигляді
\begin{equation}\label{p5-16}
\widetilde{\Phi}_1(\widetilde{\zeta}(1))=F_1(\xi_1)+
\Big((a_3y+b_3z)F_1'(\xi_1)+G_3(\xi_1)\Big)I_3\,.
\end{equation}
Помноживши рівність (\ref{p5-16}) на $I_1$\,, переконуємось у
справедливості рівності (\ref{p5-14}).

Розглянемо випадок 2). Дійсно, для трійки (\ref{har-hv==}) маємо
$$\widetilde{\zeta}(2)=x+y\widetilde{e}_2(2)+z\widetilde{e}_3(2)=x+a_2y+b_2z+a_3xI_3+b_3yI_3\,.$$
Очевидно, що
$f(\widetilde{\zeta}(2))=x+a_2y+b_2z=\xi_2\equiv\widetilde{\xi}$,
$D_2\equiv D$. Покладемо $\widetilde{F}\equiv F_2$,
$\widetilde{G}\equiv G_3$ в $D$. Тоді рівність (\ref{p5-12})
перепишеться у вигляді
\begin{equation}\label{p5-17}
\widetilde{\Phi}_2(\widetilde{\zeta}(2))=F_2(\xi_2)+
\Big((a_3y+b_3z)F_1'(\xi_2)+G_3(\xi_2)\Big)I_3\,.
\end{equation}
Помноживши рівність (\ref{p5-17}) на $I_2$\,, переконуємось у
справедливості рівності (\ref{p5-15}).

Таким чином, справедлива рівність (\ref{p5-10}):
$$
\Phi(\zeta)=I_1\,\widetilde{\Phi}_1\big(\widetilde{\zeta}(1)\big)+I_2\,\widetilde{\Phi}_2\big(\widetilde{\zeta}(2)\big),
$$
де $\Phi$ приймає значення в алгебрі $\mathbb{A}_3^2$\,, а
$\widetilde{\Phi}_1\big(\widetilde{\zeta}(1)\big),
\widetilde{\Phi}_2\big(\widetilde{\zeta}(2)\big)$ приймають значення
в $\mathbb{B}$.
\end{example}

\begin{example}
Отже, розглядаємо алгебру $\mathbb{A}_5^3$ з таблицею множення
(\ref{ex2}). Для алгебри $\mathbb{A}_5^3$ алгеброю виду $1\oplus_s
N$ є алгебра $\mathbb{A}_4$ з таблицею множення (\ref{ex21}).

Відповідно до представлення (\ref{Teor--1}), кожна моногенна функція
$\Phi$ зі значеннями в алгебрі $\mathbb{A}_5^3$ подається у вигляді
$$
\Phi(\zeta)=F_1(\xi_1)I_1+F_2(\xi_2)I_2+F_3(\xi_3)I_3+\Big((a_4y+b_4z)F_3'(\xi_3)+G_3(\xi_3)\Big)I_4+$$
\begin{equation}\label{p5-18}
+\Big((a_5y+b_5z)F_1'(\xi_1)+G_5(\xi_1)\Big)I_5
\end{equation}
$$\forall\,\zeta\in\Pi_\zeta\,,\quad \xi_u=x+a_uy+b_uz,\,\,u=1,2,3,
$$
де $F_1\,,G_5 $ --- деякі голоморфні функції в області $D_1$\,,
$F_2$ --- деяка голоморфна функція в області $D_2$\,, а $F_3\,,G_3$
--- деякі голоморфні функції в області $D_3$\,. Оскільки для
$\mathbb{A}_5^3$ \, $m=3$, то геометрично область $\Pi_\zeta$ є
перетином трьох нескінченних циліндрів:\,
$\Pi_\zeta=\Pi_\zeta(1)\cap\Pi_\zeta(2)\cap\Pi_\zeta(3)$.

Подамо функцію (\ref{p5-18}) у вигляді (\ref{p5-1}):
\begin{equation}\label{p5-19}
\Phi(\zeta)=\Phi(\zeta)I_1+\Phi(\zeta)I_2+\Phi(\zeta)I_3=:\Phi_1(\zeta)+\Phi_2(\zeta)+\Phi_3(\zeta),
\end{equation}
де
$$\Phi_1(\zeta)=F_1(\xi_1)I_1+\Big((a_5y+b_5z)F_1'(\xi_1)+G_5(\xi_1)\Big)I_5 $$
--- моногенна функція в циліндрі $\Pi_\zeta(1)$,
 функція $\Phi_2(\zeta)=F_2(\xi_2)I_2$ моногенна в циліндрі
 $\Pi_\zeta(2)$\,,
а функція
$$\Phi_3(\zeta)=F_3(\xi_3)I_3+\Big((a_4y+b_4z)F_3'(\xi_3)+G_3(\xi_3)\Big)I_4$$
моногенна в циліндрі $\Pi_\zeta(3)$.\vskip2mm

Перейдемо до розгляду моногенних функцій в алгебрі $\mathbb{A}_4$. З
представлення (\ref{Teor--1}) випливає, що кожна моногенна функція
$\widetilde{\Phi}$ зі значеннями в алгебрі $\mathbb{A}_4$ подається
у вигляді
$$
\widetilde{\Phi}(\widetilde{\zeta})=\widetilde{F}(\widetilde{\xi})+
\Big((a_4y+b_4z)\widetilde{F}'(\widetilde{\xi})+\widetilde{G}_3(\widetilde{\xi})\Big)I_4+
$$
\begin{equation}\label{p5-20}
+\Big((a_5y+b_5z)\widetilde{F}'(\widetilde{\xi})+\widetilde{G}_5(\widetilde{\xi})\Big)I_5
\quad\forall\,\widetilde{\zeta}\in\widetilde{\Pi}_{\widetilde{\zeta}}\,,\,\,\widetilde{\xi}=f(\widetilde{\zeta}),
\end{equation}
 де $\widetilde{F},\widetilde{G}_3\,,\widetilde{G}_5$
--- деякі голоморфні функції в області $D\subset\mathbb{C}$. Область $\widetilde{\Pi}_{\widetilde{\zeta}}$ є
нескінченним циліндром.

Для трійки (\ref{har-hv==+}) алгебри $\mathbb{A}_4$ маємо
$$\widetilde{\zeta}(1)=x+y\widetilde{e}_2(1)+z\widetilde{e}_3(1)=x+a_1y+b_1z+a_5xI_5+b_5yI_5\,.$$
Очевидно, що
$f(\widetilde{\zeta}(1))=x+a_1y+b_1z=\xi_1\equiv\widetilde{\xi}$,
$D_1\equiv D$. Покладемо $\widetilde{F}\equiv F_1$,
$\widetilde{G}_3\equiv G_3$, $\widetilde{G}_5\equiv G_5$ в $D$. Тоді
рівність (\ref{p5-20}) перепишеться у вигляді
$$
\widetilde{\Phi}_1(\widetilde{\zeta}(1))=F_1(\xi_1)+
\Big((a_4y+b_4z)F_1'(\xi_1)+G_3(\xi_1)\Big)I_4+
$$
\begin{equation}\label{p5-21}
+\Big((a_5y+b_5z)F_1'(\xi_1)+G_5(\xi_1)\Big)I_5\,.
\end{equation}
Помноживши рівність (\ref{p5-21}) на $I_1$\,, переконуємось у
справедливості рівності
$$
I_1\,\widetilde{\Phi}_1\big(\widetilde{\zeta}(1)\big)=\Phi_1(\zeta)\,.
$$

Для трійки (\ref{har-hv==++}) алгебри $\mathbb{A}_4$  маємо
$$\widetilde{\zeta}(2)=x+y\widetilde{e}_2(2)+z\widetilde{e}_3(2)=x+a_2y+b_2z.$$
Очевидно, що
$f(\widetilde{\zeta}(2))=\widetilde{\zeta}(2)=x+a_2y+b_2z=\xi_2\equiv\widetilde{\xi}$,
$D_2\equiv D$. Покладемо $\widetilde{F}\equiv F_2$,
$\widetilde{G}_3\equiv G_3$, $\widetilde{G}_5\equiv G_5$ в $D$. Тоді
рівність (\ref{p5-20}) набуде вигляду
$$
\widetilde{\Phi}_2(\widetilde{\zeta}(2))=F_2(\xi_2)+
\Big((a_4y+b_4z)F_2'(\xi_2)+G_3(\xi_2)\Big)I_4+
$$
\begin{equation}\label{p5-22}
+\Big((a_5y+b_5z)F_2'(\xi_2)+G_5(\xi_2)\Big)I_5\,.
\end{equation}
Помноживши рівність (\ref{p5-22}) на $I_2$\,, переконуємось у
справедливості рівності
$$
I_2\,\widetilde{\Phi}_2\big(\widetilde{\zeta}(2)\big)=\Phi_2(\zeta)\,.
$$

Нарешті, для трійки (\ref{har-hv==+++}) отримуємо
$$\widetilde{\zeta}(3)=x+y\widetilde{e}_2(3)+z\widetilde{e}_3(3)=x+a_3y+b_3z+a_4xI_4+b_4yI_4\,.$$
Очевидно, що
$f(\widetilde{\zeta}(3))=x+a_3y+b_3z=\xi_3\equiv\widetilde{\xi}$,
$D_3\equiv D$. Покладемо $\widetilde{F}\equiv F_3$,
$\widetilde{G}_3\equiv G_3$, $\widetilde{G}_5\equiv G_5$ в $D$. Тоді
рівність (\ref{p5-20}) перепишеться у вигляді
$$
\widetilde{\Phi}_3(\widetilde{\zeta}(3))=F_3(\xi_3)+
\Big((a_4y+b_4z)F_3'(\xi_3)+G_3(\xi_3)\Big)I_4+
$$
\begin{equation}\label{p5-23}
+\Big((a_5y+b_5z)F_3'(\xi_3)+G_5(\xi_3)\Big)I_5\,.
\end{equation}
Помноживши рівність (\ref{p5-23}) на $I_3$\,, переконуємось у
справедливості рівності
$$
I_3\,\widetilde{\Phi}_3\big(\widetilde{\zeta}(3)\big)=\Phi_3(\zeta)\,.
$$

Таким чином, справедлива рівність (\ref{p5-10}):
$$
\Phi(\zeta)=I_1\,\widetilde{\Phi}_1\big(\widetilde{\zeta}(1)\big)+I_2\,\widetilde{\Phi}_2\big(\widetilde{\zeta}(2)\big)+
I_3\,\widetilde{\Phi}_3\big(\widetilde{\zeta}(3)\big),
$$
де $\Phi$ приймає значення в алгебрі $\mathbb{A}_5^3$\,, а
$\widetilde{\Phi}_1\big(\widetilde{\zeta}(1)\big)$,
$\widetilde{\Phi}_2\big(\widetilde{\zeta}(2)\big)$,
$\widetilde{\Phi}_3\big(\widetilde{\zeta}(3)\big)$  приймають
значення в $\mathbb{A}_4$.
\end{example}

\bigskip

ІНФОРМАЦІЯ ПРО АВТОРА

\medskip
Віталій Станіславович Шпаківський\\
Інститут математики НАН України,\\
вул. Терещенківська, 3, Київ, Україна\\
shpakivskyi86@gmail.com
\end{document}